\documentclass[a4paper,11pt]{article}
\usepackage{amsmath, amsfonts, amscd, amssymb, amsthm, enumerate, xypic}

\newcommand{\Mat}{\operatorname{M}}

\newcommand{\Mata}{\operatorname{A}}

\newcommand{\GL}{\operatorname{GL}}

\newcommand{\Ker}{\operatorname{Ker}}

\newcommand{\urk}{\operatorname{urk}}

\def\K{\mathbb{K}}

\def\calA{\mathcal{A}}

\def\calM{\mathcal{M}}

\def\calS{\mathcal{S}}

\def\calV{\mathcal{V}}
\def\calW{\mathcal{W}}

\def\lcro{\mathopen{[\![}}
\def\rcro{\mathclose{]\!]}}

\theoremstyle{definition}

\theoremstyle{plain}
\newtheorem{theo}{Theorem}

\theoremstyle{plain}

\theoremstyle{remark}

\title{From primitive spaces of bounded rank matrices to a generalized Gerstenhaber theorem}
\author{Cl\'ement de Seguins Pazzis\footnote{Universit\'e de Versailles Saint-Quentin-en-Yvelines, Laboratoire de Math\'ematiques
de Versailles, 45 avenue des Etats-Unis, 78035 Versailles cedex, France}
\footnote{e-mail address: dsp.prof@gmail.com}}

\begin{document}

\thispagestyle{plain}

\maketitle

\begin{abstract}
A recent generalization of Gerstenhaber's theorem on spaces of nilpotent matrices is derived, under
mild conditions on the cardinality of the underlying field, from Atkinson's structure theorem on
primitive spaces of bounded rank matrices.
\end{abstract}



\section{Introduction}

In the modern geometric theory of matrices, two deep structure theorems stand out:
Dieudonn\'e's theorem \cite{Dieudonne} on large spaces of singular matrices - later generalized
by Atkinson and Lloyd \cite{AtkLloydLarge} - and Gerstenhaber's theorem on spaces of nilpotent matrices \cite{Gerstenhaber}.
Very recent advances have been made in both problems: Atkinson and Lloyd's extension
of Dieudonn\'e's theorem has been shown to hold for almost all fields \cite{dSPclass},
while there have been several generalizations of Gerstenhaber's theorem, most notably
to \textbf{trivial spectrum} spaces of matrices, i.e.,
subspaces of square matrices having no non-zero eigenvalue in their field of definition.

In our view, another structure theorem has not received the attention it deserved from the mathematics community:
it is Atkinson's theorem on primitive spaces of bounded rank matrices \cite[Theorem B]{AtkinsonPrim}. In this note, we will expose
a yet unknown connection between Atkinson's theorem and Gerstenhaber's, and will
use this insight to obtain a greatly simplified proof of
a slightly weaker version of the classification theorem for trivial spectrum spaces of matrices.
We believe that the technique displayed here has a great potential to deliver new insights into the structure of
spaces of nilpotent matrices, as it only uses a limit case in Atkinson's theorem.

\section{Notation and basic definition}

Here, $\K$ denotes an arbitrary field, and $\Mat_{m,n}(\K)$, $\Mat_n(\K)$, $\Mata_n(\K)$, $\GL_n(\K)$
denote, respectively, the sets of all $m \times n$ matrices, $n \times n$ matrices, $n \times n$ alternating matrices,
and $n \times n$ invertible matrices with entries in $\K$. Two subsets $\calV$ and $\calW$ of $\Mat_{m,n}(\K)$
are called equivalent when there exists a pair $(P,Q)\in \GL_m(\K) \times \GL_n(\K)$ such that $\calV=P\,\calW\,Q$,
which means that $\calV$ and $\calW$ represent the same set of linear operators in different choices of bases of the source and goal spaces.
If $m=n$, we say that $\calV$ and $\calW$ are similar when, in the above condition, we require that $Q=P^{-1}$,
meaning that $\calV$ and $\calW$ represent the same set of endomorphisms of a vector space in two potentially different bases.

Given a subset $\calV$ of $\Mat_{m,n}(\K)$, the upper rank of $\calV$, denoted by $\urk(\calV)$,
is the largest rank for a matrix in $\calV$. Note that two equivalent subsets share the same upper rank.

A linear subspace $\calV$ of $\Mat_{m,n}(\K)$ is called \textbf{primitive} when it satisfies the following four conditions:
\begin{enumerate}[(i)]
\item $\calV$ is not equivalent to a space of matrices with the last column equal to zero;
\item $\calV$ is not equivalent to a space of matrices with the last row equal to zero;
\item $\calV$ is not equivalent to a space $\calV'$ in which every matrix is written as
$M=\begin{bmatrix}
H(M) & [?]_{m \times 1}
\end{bmatrix}$ and $\urk H(\calV')<\urk \calV$;
\item $\calV$ is not equivalent to a space $\calV'$ in which every matrix is written as
$M=\begin{bmatrix}
R(M) \\
[?]_{1 \times n}
\end{bmatrix}$ and $\urk R(\calV')<\urk \calV$.
\end{enumerate}
Moreover, we say that $\calV$ is \textbf{semi-primitive} when it is only required to satisfy conditions (i), (ii) and (iii);
we say that $\calV$ is \textbf{reduced} when it is only required to satisfy conditions (i) and (ii).

We note that the upper rank of a semi-primitive subspace of $\Mat_{m,n}(\K)$ is always less than $n$,
so that the upper rank of a primitive subspace is always less than $m$ and $n$.
As shown by the first statement in Theorem 1 of \cite{AtkLloydPrim}, the primitive spaces are the elementary pieces
upon which are built all those matrix spaces with upper rank less than the number of rows and the number of columns.
A fundamental example of primitive space can be derived from
the canonical pairing
$$\varphi_n : \K^n \times \K^n \rightarrow \K^n \wedge \K^n.$$
In the canonical basis $(e_1,\dots,e_n)$ of $\K^n$ and the lexicographically ordered basis $(e_i \wedge e_j)_{1 \leq i<j \leq n}$
of $\K^n \wedge \K^n$, one then takes the space $S(\varphi_n)$
of all matrices representing the operators $x \wedge -$ for $x \in \K^n$.
One checks that, for $n \geq 2$, $S(\varphi_n)$ is a semi-primitive linear subspace of $\Mat_{\binom{n}{2},n}(\K)$ -
and even a primitive one if $n>2$ - and that it is also the space of all
matrices
$$\begin{bmatrix}
X^T A_1 \\
\vdots \\
X^T A_{\binom{n}{2}}
\end{bmatrix} \quad \text{with $X \in \K^n$,}$$
where the matrices $A_1,\dots,A_{\binom{n}{2}}$ are the elements of $(e_i e_j^T-e_j e_i^T)_{1 \leq i<j \leq n}$ put in the lexicographical order.
A subspace $\calV$ equivalent to $S(\varphi_n)$ is exactly a subspace for which there
is a basis $(B_1,\dots,B_{\binom{n}{2}})$ of $\Mata_n(\K)$ together with some $P \in \GL_n(\K)$ such that
$\calV$ is the space of all matrices
$$\begin{bmatrix}
X^T B_1P \\
\vdots \\
X^T B_{\binom{n}{2}}P
\end{bmatrix} \quad \text{with $X \in \K^n$.}$$

\newpage
\section{Atkinson's theorem}

We now state a special case of the transposed version of Atkinson's theorem, combining Lemma 6 of \cite{AtkLloydPrim}
and Theorem B of \cite{AtkinsonPrim}.

\begin{theo}[Atkinson]\label{Atkinson}
Let $\calS$ be a semi-primitive linear subspace of $\Mat_{m,n}(\K)$. Set $r:=\urk \calS$ and assume that $\# \K>r$.
Then, $m \leq \frac{r(r+1)}{2}$. If in addition $m=\frac{r(r+1)}{2}$ and $r>1$, then $n=r+1$ and
$\calS$ is equivalent to $S(\varphi_n)$.
\end{theo}

Note that, for the case $m=\frac{r(r+1)}{2}$ and $r>1$, Atkinson only states his result for primitive spaces:
however, if we put Atkinson's arguments in the context of our version (i.e., we transpose them), then
the only instance when he uses condition (iv) is to discard the case where $\calS$ might be
equivalent to a space $\calS'$ of matrices of the form
$M=\begin{bmatrix}
[?]_{1 \times q} & [?]_{1 \times (n-q)} \\
H(M) & [0]_{(m-1) \times (n-q)}
\end{bmatrix}$
with $q \in \lcro 0,n-1\rcro$ and a reduced space $H(\calS')$. But then
$\urk H(\calS')\leq r-1$ and Lemma 6 of \cite{AtkLloydPrim} would show that $m-1 \leq \binom{r}{2}<\binom{r+1}{2}-1$ as,
in the terminology of Atkinson and Lloyd, the space $H(\calS')$ would have the column property.

\vskip 2mm
We also remark that the conclusion of the second statement still holds in the case $m=1$ and $n=2$, as then
$\calS=\Mat_{1,2}(\K)=S(\varphi_2)$.

\section{The generalized Gerstenhaber theorem}

Seemingly unrelated to Atkinson's theory of primitive spaces is the following generalization
of Gerstenhaber's theorem, proved in \cite{dSPlargeaffine},
in which a linear subspace $\calV$ of $\Mat_n(\K)$ is called \textbf{irreducible} if
no proper non-zero linear subspace of $\K^n$ is globally invariant under all the elements of~$\calV$.

\begin{theo}[Generalized Gerstenhaber theorem]\label{supergerstenhaber}
Let $\calV$ be a trivial spectrum linear subspace of $\Mat_n(\K)$.
Then,
\begin{equation}\label{majoration}
\dim \calV \leq \frac{n(n-1)}{2}\cdot
\end{equation}
If equality holds in \eqref{majoration} with $\# \K \geq 3$ and
$\calV$ irreducible, then $\calV=P\Mata_n(\K)$ for some $P \in \GL_n(\K)$.
\end{theo}

Given $P \in \GL_n(\K)$, one then shows that $P\Mata_n(\K)$ is an irreducible space
with a trivial spectrum if and only if $X \mapsto X^TPX$ is a non-isotropic quadratic form on $\K^n$,
see \cite[Lemma 10]{dSPlargeaffine}. Theorem \ref{supergerstenhaber} may be used to provide a complete classification of trivial spectrum spaces
with the maximal dimension - up to similarity - reducing the problem to the classification of non-isotropic bilinear forms
\cite[Theorem 4]{dSPlargeaffine} up to similarity.

The main rationale for studying trivial spectrum spaces
stems from their ties to affine spaces of non-singular matrices.
In short, if $\calA$ is an affine subspace of non-singular matrices of $\Mat_n(\K)$ which contains $I_n$,
then its translation vector space $\calA-I_n$ is a trivial spectrum linear subspace.
With little effort, Theorem \ref{supergerstenhaber} shows that any such affine subspace has dimension less than or equal to
$\binom{n}{2}$, and reduces the classification - up to equivalence - of those with the maximal dimension to
the classification of non-isotropic quadratic forms on $\K$, up to similarity \cite[Theorem 7]{dSPlargeaffine}.
This has led to similar structure theorems for large affine spaces of matrices with a lower bound on the rank
\cite{dSPlargeaffinerank}.

Finally, Theorem \ref{supergerstenhaber} yields Gerstenhaber's theorem - which we recall now - as an easy consequence
for fields with more than two elements (for proofs that hold regardless of the cardinality of the field,
see \cite{dSPGerstenhaber} and \cite{Serezhkin}, while an elegant proof that works only for fields with more than two elements
is featured in \cite{Mathes}).

\begin{theo}[Gerstenhaber's theorem \cite{Gerstenhaber}]
Let $\calV$ be a linear subspace of $\Mat_n(\K)$ in which all the matrices are nilpotent.
Then, $\dim \calV \leq \frac{n(n-1)}{2}$, and equality holds if and only if $\calV$
is similar to the space of all strictly upper-triangular $n\times n$ matrices.
\end{theo}

\section{Proof of the generalized Gerstenhaber theorem in a restricted setting}

Spread over 20 dense pages, the only known proof of Theorem \ref{supergerstenhaber}
is a very intricate tour de force. Using Atkinson's theorem on semi-primitive spaces, we shall now
give a much shorter alternative proof \emph{under the additional assumption that $\# \K \geq n$}.

Let us immediately explain the connection with Atkinson's theorem.
Let $\calV$ be a trivial spectrum subspace of $\Mat_n(\K)$.
Then, $\calV$ has an interesting property that was the cornerstone of the proof of Theorem \ref{supergerstenhaber}
featured in \cite{dSPlargeaffine}: it is \textbf{totally intransitive} in the sense that, for every non-zero vector $X \in \K^n$,
the linear subspace $\calV X:=\{NX \mid N \in \calV\}$ is a proper subspace of $\K^n$ as it cannot contain $X$.

We now combine this simple fact with a duality argument to create a linear subspace of $\Mat_{m,n}(\K)$ with upper rank less than $n$, where
$m:=\dim \calV$. For $X \in \K^n$, consider the bilinear form
$$\widehat{X}: \; (N,Y) \in \calV \times \K^n \, \longmapsto \,  Y^TNX \in \K,$$
and denote by $\calM \subset \Mat_{m,n}(\K)$ the space of all matrices representing such forms
in chosen bases of $\calV$ and $\K^n$.
In other words, if the chosen basis of $\calV$ is $(B_1,\dots,B_m)$,
then there is an invertible matrix $P \in \GL_n(\K)$ such that, for all $X \in \K^n$, the matrix of $\widehat{X}$ is
$$\begin{bmatrix}
X^T B_1^TP \\
\vdots \\
X^T B_m^TP
\end{bmatrix}.$$
Note that $\calV$ always satisfies condition (ii), as
a matrix $N \in \Mat_n(\K)$ which satisfies $Y^T NX=0$ for all $(X,Y)\in (\K^n)^2$ is necessarily zero.
We also see that $\urk \calV <n$ since, given a non-zero vector $X \in \K^n$,
the set $\calV X$ is a proper linear subspace of $\K^n$, which yields a non-zero vector $Y \in \K^n$ for which $Y^T\calV X=0$.

Now, it may very well happen that $\calM$ is not semi-primitive, but let us assume for the moment that it is.
Then, by Atkinson's theorem, we would have $m \leq \binom{n}{2}$, and, in case of equality, we would
find a basis $(A_1,\dots,A_m)$ of $\Mata_n(\K)$ together with an invertible matrix $Q \in \GL_n(\K)$ such that
$$\forall X \in \K^n, \; \begin{bmatrix}
X^T B_1^TP \\
\vdots \\
X^T B_m^TP
\end{bmatrix}=\begin{bmatrix}
X^T A_1Q \\
\vdots \\
X^T A_mQ
\end{bmatrix};$$
this would yield $B_i^TP=A_iQ$ for all $i \in \lcro 1,m\rcro$, and hence we would conclude that
$$\calV=(QP^{-1})^T\Mata_n(\K).$$

\paragraph{}
Let us now return to the general case. The proof works by induction on $n$.
The result is trivial for $n=1$ and we now assume that $n \geq 2$.
Note first that we may always assume that $\calV$ is irreducible, for if it is not,
then, replacing $\calV$ by a similar subspace,
we can see that no generality is lost in assuming that there exists an integer $p \in \lcro 1,n-1\rcro$ such that
every matrix of $\calV$ splits up as
$$N=\begin{bmatrix}
A(N) & [?]_{p \times (n-p)} \\
[0] & B(N)
\end{bmatrix},$$
where $A(N)$ and $B(N)$ are, respectively, $p \times p$ and $(n-p) \times (n-p)$ matrices;
$A(\calV)$ and $B(\calV)$ are then trivial spectrum subspaces, respectively, of $\Mat_p(\K)$ and $\Mat_{n-p}(\K)$, and by induction we deduce that
$$\dim \calV \leq \dbinom{p}{2}+\dbinom{n-p}{2}+p(n-p)=\dbinom{n}{2}.$$
In the rest of the proof, we assume that $\calV$ is irreducible.

Let us come back to the matrix space $\calM$. It now satisfies condition (i), for if it did not,
then we would have a non-zero vector $Y \in \K^n$ for which $\forall X \in \K^n, \; \forall N \in \calV, \; Y^T NX=0$,
and hence all the elements of $\calV$ would map all the vectors of $\K^n$ into some fixed
linear hyperplane of $\K^n$, contradicting the assumed irreducibility of $\calV$.
Therefore, $\calM$ satisfies conditions (i), (ii), and $\urk \calM<n$. If $\calM$ is semi-primitive, then
the conclusion follows as explained earlier.

\paragraph{}If we now assume that $\calM$ is not semi-primitive, then
we may find a \emph{minimal} integer $d \in \lcro 1,n-1\rcro$ for which the bases of $\calV$
and $\K^n$ are chosen so that every matrix of $\calM$ has the form
$$M=\begin{bmatrix}
H(M) & [?]_{m \times (n-d)}
\end{bmatrix},$$
and $H(\calM)$ is a linear subspace of $\Mat_{m,d}(\K)$ with upper rank less than $d$.
In this situation, we may also modify the basis of $\calV$ further to the point where, for every
$M \in \calM$, we have
$$H(M)=\begin{bmatrix}
K(M) \\
[0]_{(m-c) \times d}
\end{bmatrix},$$
where $K(\calM)$ is a linear subspace of $\Mat_{c,d}(\K)$ with upper rank less than $d$
and which satisfies condition (ii). Then, using the minimality of $d$ and the fact that $\calM$ satisfies condition (i),
we can see that $K(\calM)$ is a semi-primitive subspace of $\Mat_{c,d}(\K)$.
Theorem \ref{Atkinson} then yields
$$c \leq \binom{d}{2}.$$

Let us now come back to $\calV$ and see how the above reduction plays out.
As we may replace $\calV$ with a similar subspace, no generality is lost in assuming that
the chosen basis of $\K^n$ is the canonical one.
Let us then write every matrix $N$ of $\calV$ as
$$N=\begin{bmatrix}
R_1(N) \\
R_2(N)
\end{bmatrix}=\begin{bmatrix}
A(N) & C(N) \\
B(N) & D(N)
\end{bmatrix},$$
where $R_1(N)$, $R_2(N)$, $A(N)$, $B(N)$, $C(N)$ and $D(N)$ are, respectively,
$d \times n$, $(n-d) \times n$, $d \times d$, $(n-d) \times d$, $d \times (n-d)$ and $(n-d) \times (n-d)$ matrices.
Set also
$$\calW:=\Ker R_1.$$
As $K(\calM)$ satisfies condition (ii), one sees that the last $m-c$ vectors of the chosen basis of $\calV$
span the subspace of all matrices of $\calV$ in which the first $d$ rows equal zero, and hence
$$m-c=\dim \calW.$$
On the other hand, every matrix $N$ of $\calW$ splits up as
$$N=\begin{bmatrix}
[0] & [0] \\
B(N) & D(N)
\end{bmatrix},$$
and hence $D(\calW)$ is a trivial spectrum subspace of $\Mat_{n-d}(\K)$.
By induction, we have $\dim D(\calW) \leq \dbinom{n-d}{2}$, and therefore
$$\dim \calW \leq (n-d)\,d+\binom{n-d}{2}.$$
We conclude that
$$m = c+\dim \calW \leq \binom{d}{2}+(n-d)\,d+\binom{n-d}{2}=\binom{n}{2},$$
thus completing the proof of inequality \eqref{majoration}.

If we now assume that $m=\dbinom{n}{2}$ on top of the previous assumptions, then all the above inequalities
turn out to be equalities, and in particular we have:
\begin{enumerate}[(a)]
\item $\dim D(\calW)=\dbinom{n-d}{2}$;
\item The space $\calV$ contains every matrix of the form
$$\begin{bmatrix}
[0]_{d \times d} & [0]_{d \times (n-d)} \\
[?]_{(n-d) \times d} & [0]_{(n-d) \times (n-d)}
\end{bmatrix}.$$
\end{enumerate}
Using point (b), we deduce that, for every $N \in \calV$, the matrix
$$\begin{bmatrix}
A(N) & C(N) \\
[0] & D(N)
\end{bmatrix}$$
belongs to $\calV$, and hence $D(\calV)$ is a trivial spectrum subspace of $\Mat_{n-d}(\K)$.
By induction, we have $\dim D(\calV)\leq \binom{n-d}{2}=\dim D(\calW)$ with $D(\calW) \subset D(\calV)$, which yields
$D(\calW)=D(\calV)$.

We shall now obtain a contradiction from an invariance argument.
Assume that some $N_0 \in \calV$ satisfies $C(N_0) \neq 0$.
Then, we may find a non-zero vector $x \in \K^{n-d}$ such that $C(N_0)\,x \neq 0$, and then choose
$R \in \Mat_{n-d,d}(\K)$ for which $R\,C(N_0)\,x=x$.
With the invertible matrix
$P:=\begin{bmatrix}
I_d & [0] \\
R & I_{n-d}
\end{bmatrix}$, one computes that
$$\forall N \in \calV, \quad PNP^{-1}=\begin{bmatrix}
A(N)-C(N)R & C(N) \\
B(N)+RA(N)-RC(N)R & D(N)+R C(N)
\end{bmatrix}.$$
In the new trivial spectrum space $\calV':=P\calV P^{-1}$, the subspace of matrices with all first $d$ rows equal to zero
is still $\calW$.
Thus, with the above dimensional arguments applied to $\calV'$, we can deduce that $D(N_0)+R\,C(N_0)$ belongs to $D(\calV')=D(\calW)=D(\calV)$,
and hence $R\,C(N_0)$ belongs to $D(\calV)$. This is absurd because we have seen that $D(\calV)$ is a trivial spectrum subspace of $\Mat_{n-d}(\K)$.

Therefore, $C(N)=0$ for all $N \in \calV$, and hence $\{0\} \times \K^{n-d}$ is
a globally invariant subspace for all the matrices of $\calV$.
This contradicts the assumed irreducibility of $\calV$ and concludes the proof.

\end{document}